\documentclass[12pt,reqno]{amsart}

\usepackage{amssymb,amsbsy,amsthm,amsmath,graphicx,epsfig}

\usepackage{amssymb}
\usepackage{amscd}
\usepackage{amsfonts}
\usepackage{mathrsfs}
\usepackage{setspace}
\usepackage{version}

\theoremstyle{plain}
\numberwithin{equation}{subsection} \theoremstyle{plain}
\newtheorem{theorem}[subsection]{Theorem}

\newtheorem{example}[subsection]{Example}

\newtheorem{remark}[subsection]{Remark}

\newtheorem{rmm}[subsection]{}

\renewcommand{\leq}{\leqslant}
\renewcommand{\geq}{\geqslant}
\newsavebox{\proofbox}
\savebox{\proofbox}{\begin{picture}(7,7)  \put(0,0){\framebox(7,7){}}\end{picture}}

\newcommand\Z{\mathbb{Z}}
\newcommand\R{\mathbb{R}}

\newcommand\C{\mathbb{C}}

\newcommand\HH{\mathbb{H}}

\newcommand\SL{\operatorname{SL}}

\newcommand\Q{\mathbb{Q}}

\renewcommand\P{\mathcal{P}}

\begin{document}

\title[Salem numbers and $\lambda_1$]{Salem numbers and the spectrum of hyperbolic surfaces}

\author{Emmanuel Breuillard and Bertrand Deroin}

\keywords{Cheeger constant, first eigenvalue of the Laplacian, Salem number, arithmetic surfaces}

\address{Mathematisches Institut\\
WWU Muenster\\
62 Einsteinstrasse\\
48149 Muenster\\
GERMANY}
\email{emmanuel.breuillard@uni-muenster.de}

\address{CNRS\\
AGM\\
2 av. Adolphe Chauvin\\
95302 Cergy Pontoise Cedex\\
FRANCE}
\email{bertrand.deroin@u-cergy.fr}

\begin{abstract} We give a reformulation of Salem's conjecture about the absence of Salem numbers near one in terms of a uniform spectral gap for certain arithmetic hyperbolic surfaces.
\end{abstract}

\date{December 2016}
\maketitle

\bigskip

\bigskip

\section{introduction}
A Salem number is an algebraic integer with only one Galois conjugate inside the open unit disc and at least one on the unit circle. Salem numbers are real and in $(1,+\infty)$. Salem conjectured \cite{salem-duke, salem-book} that there is an absolute constant $\tau>1$ such that every Salem number is $\geq \tau$. We call this the ``no small Salem number conjecture'', or \emph{Salem's conjecture} for short.

It was observed some time ago by Sury \cite{sury} that the Salem conjecture has the following beautiful geometric reformulation:

\begin{rmm}\label{salemeq} \noindent {\it Salem's conjecture holds if and only if there is a uniform positive lower bound on the length of closed geodesics in arithmetic hyperbolic $2$-orbifolds.}
\end{rmm}


We also refer the reader to \cite[Chapter IX (4.21)]{margulis}, \cite{ghate-hironaka},\cite[Chapter 12]{maclachlan-reid} and to Section \ref{salem-sec} below for an exposition of this phenomenon. By a hyperbolic $2$-orbifold, we mean the quotient of the hyperbolic plane $\HH^2$ by a discrete subgroup $\Gamma$ of $PSL_2(\R)$ of finite co-volume. If $\Gamma$ is torsion-free, we will call the corresponding quotient $\HH^2/\Gamma$ a hyperbolic surface. It is arithmetic iff $\Gamma$ is commensurable to a congruence lattice, that is a discrete subgroup of $PSL_2(\R)$ which is a congruence subgroup of the group of units of an order in a quaternion algebra defined over a totally real number field (see e.g. \cite{takeuchi}).

The length of a closed geodesic is the translation length $L_\gamma$ of a hyperbolic element $\gamma \in \Gamma$ on its axis. For arithmetic hyperbolic $2$-orbifolds, $\exp{L_\gamma}$ is always a Salem number. In particular, if Salem's conjecture holds, then the uniform lower bound on the length of closed geodesics holds for all arithmetic hyperbolic $2$-orbifolds. Conversely, every Salem number $\alpha$ arises this way. Namely $\alpha$ is a quadratic unit over a certain totally real field, which can be used to build a quaternion algebra giving rise to an arithmetic group admitting a hyperbolic element of translation length $2\log \alpha$ (see \cite{sury,ghate-hironaka}).  We refer the reader to \cite{maclachlan-reid} for background on arithmetic surfaces.

As it turns out a uniform lower bound on the length of closed geodesics in congruence arithmetic hyperbolic surfaces (i.e. with $\Gamma$ congruence and torsion-free) is enough to yield Salem's conjecture (see (\ref{salem2smoothsurface}) below).

\bigskip

Congruence $2$-orbifolds are well-known to satisfy the universal lower bound $\lambda_1 \geq \frac{3}{16}$ (this follows from the work of Gelbart-Jacquet \cite{gelbart-jacquet} and Jacquet-Langlands \cite{jacquet-langlands} as observed by Vigneras in \cite{vigneras}). The goal of this note is to show how this information can be combined with Sury's geometric criterion $(\ref{salemeq})$ in order to obtain the following spectral reformulation of Salem's conjecture:

\begin{theorem}\label{specriterion} \noindent {\it Salem's conjecture holds if and only if there is a uniform $c>0$ such that
\begin{equation}\label{ineqqqq}\lambda_1(\widetilde{\Sigma}) \geq \frac{c}{area(\widetilde{\Sigma})}\end{equation} for all $2$-covers $\widetilde{\Sigma}$ of all compact congruence arithmetic hyperbolic $2$-orbifolds $\Sigma$.}
\end{theorem}

We consider only unramified covers. Here again, in order to get Salem's conjecture, it is enough to know the above lower bound either when $\Sigma$ varies in the class of congruence arithmetic hyperbolic surfaces, or in the class of $2$-orbifolds of the form $\HH^2/\Gamma$, where $\Gamma$ is a maximal arithmetic subgroup of $PSL_2(\R)$. It is also enough to know it for $2$-orbifolds $\HH^2/\Gamma$, where the Fuchsian group $\Gamma$ is the subgroup $\Gamma^{1}_{\mathcal{O}}$ of norm $1$ elements in a maximal order $\mathcal{O}$ in a quaternion algebra defined over a totally real number field. In particular, we do not need to consider towers of congruence subgroups of a given maximal arithmetic group, only the base surface is needed. What is important is the uniformity over the number fields.

The above lower bound for $2$-covers of non compact congruence surfaces is true and easy to establish, because non co-compact arithmetic subgroups are commensurable to $PSL_2(\Z)$.

\bigskip

The link between Sury's $(\ref{salemeq})$ and Theorem \ref{specriterion} is provided by the Cheeger-Buser inequality. Recall that if
 $M$ is a compact connected $n$-dimensional smooth Riemannian manifold and $\lambda_1(M)$ is the first non-zero eigenvalue  of the Laplace-Beltrami operator, then the Cheeger constant $h(M)$  of the manifold is defined by:
\begin{equation}\label{cheegdef}
h(M):=\inf_A \frac{vol(\partial A)}{\min\{vol(A),vol(M\setminus A)\}},
\end{equation}
where $A$ denotes an open submanifold of $M$ with volume $vol(A)$ and a smooth boundary with $(n-1)$-dimensional area $vol_{n-1}(\partial A)$. Recall the Cheeger-Buser inequality, which holds  in presence of a Ricci curvature lower bound on $M$, namely $Ric(M) \geq -(n-1)$, where $n=\dim M$:

\begin{equation}\label{buser}
\frac{1}{4}h(M)^2 \leq \lambda_1(M) \leq  2(n-1)h(M) + 10h(M)^2.
\end{equation}

The left hand side is the Cheeger inequality \cite{cheeger}, while the right hand side is the Buser inequality \cite{buser}. Combined together these inequalities say that $\lambda_1(M)$ and $h(M)^2$ are comparable up to multiplicative constants, provided $h(M)$ is not too small. However if $h(M)$ is small, then the comparison leaves open the possibility for $\lambda_1(M)$ to be either close to $h(M)^2$, or rather close to $h(M)$. Both cases occur: for a flat circle of radius $R \gg 1$ for example, $\lambda_1(M)=\frac{4\pi}{R^2}$ and $h(M)=\frac{4}{R}$, while for a compact hyperbolic surface of genus $2$ say, $\lambda_1(M)$ is comparable to $h(M)$ up to  multiplicative constants (see e.g. \cite{schoen-wolpert-yau} or \cite{dodziuk-randol-sullivan}).

In particular Theorem \ref{specriterion} is not a direct consequence of $(\ref{salemeq})$ and the Cheeger-Buser inequality. In order to derive the equivalence in Theorem \ref{specriterion}, we show that when looking at a $2$-cover $M$ of a base manifold, and in presence of a spectral gap for the base manifold, then a strengthened version of Cheeger's inequality holds, implying that $\lambda_1(M)$ and $h(M)$ are comparable up to multiplicative constants, unless they are very large.

\begin{theorem}\label{main} Let $M$ be a compact connected Riemannian manifold and $M'$ a $2$-cover of $M$. Then
$$\lambda_1(M') \geq \frac{1}{4} \sqrt{\lambda_1(M)} \cdot h(M').$$
\end{theorem}

The proof of Theorem \ref{main} is given in Section \ref{2cover-sec}. In a follow-up paper \cite{breuillard-deroin2} we extend Theorem \ref{main} to arbitrary finite covers. Theorem \ref{specriterion} is proven in Section \ref{salem-sec}. In the remainder of this introduction, we further discuss Theorem \ref{specriterion} and its possible extension to hyperbolic $n$-manifolds for $n \geq 3$.
\bigskip

It is easy to build compact hyperbolic surfaces that violate the inequality $(\ref{ineqqqq})$, for example by pinching a simple closed geodesic, or by taking large cyclic covers of a fixed surface. On the other hand, Buser and Sarnak showed in \cite{buser-sarnak} that every compact hyperbolic surface $\Sigma$ admits a simple closed non separating geodesic of length at most $L_g=2\log(4g-2)$, where $g=g(\Sigma)$ is the genus of $\Sigma$, and thus cutting along that geodesic and gluing two copies of the surface along the cuts, one obtains a $2$-cover $\widetilde{\Sigma}$ whose Cheeger constant is at most $L_g/area(\Sigma)$, and by Buser's inequality $(\ref{buser})$:
$$\lambda_1(\widetilde{\Sigma}) \leq 8\frac{\log (area(\widetilde{\Sigma}))}{area(\widetilde{\Sigma})}.$$

\bigskip

It is worth asking what happens for arithmetic Kleinian groups, i.e. arithmetic lattices in $PSL_2(\C)$. In this case the existence of a positive lower bound on the length of closed geodesics is again equivalent to the following refinement of Salem's conjecture: if an algebraic integer has  either only one conjugate outside of the closed unit disc or has exactly two non real conjugates, then these conjugates are bounded away from $1$ in modulus by a universal positive constant, see \cite[Chapter 12]{maclachlan-reid}. However there seems to be no spectral criterion immediately analogous to Theorem $(\ref{specriterion})$ in this case. In fact, as we will show at the end of this note, the analogous lower bound on $\lambda_1$ does hold:


\begin{rmm}\label{kleineq} \noindent {\it There is a uniform $c>0$ such that
$$\lambda_1(\widetilde{\Sigma}) \geq \frac{c}{volume(\widetilde{\Sigma})}$$ for all $2$-covers $\widetilde{\Sigma}$ of all compact congruence hyperbolic $n$-manifolds for all $n \geq 3$.}
\end{rmm}

Therefore it is only in dimension $2$, where the length of closed geodesics is intimately related to the $\lambda_1$ that we get a spectral reformulation of the number theoretic conjecture.

Note once again that taking large cyclic covers of a given compact hyperbolic $n$-manifold (say with positive first Betti number), it is easy to construct examples of compact hyperbolic $n$-manifolds that violate the above lower bound on $\lambda_1$.

On the other hand, by a similar cutting and gluing procedure, it is easy to construct examples of $2$-covers of congruence hyperbolic $n$-manifolds (for each $n \geq 2$) with arbitrarily large volume, whose $\lambda_1$ also satisfies an upper bound of same order of magnitude as the lower bound in  $(\ref{kleineq})$, i.e. of order $1/vol$ up to multiplicative constants; e.g. start with a standard congruence arithmetic hyperbolic $n$-manifold given by a quadratic form on $n+1$ variables (\cite[Exemple 1, p.120]{bergeron}), pick a totally geodesic hypersurface obtained by ``forgetting a variable'' and take a suitable finite congruence cover to increase the volume if necessary while preserving the hypersurface using subgroup separability (see \cite[Lemme principal p.113]{bergeron}). We are indebted to N. Bergeron for this observation.

\bigskip


\bigskip

\noindent \emph{Acknowledgements.} We are grateful to G. Chenevier, N. Bergeron, A. Reid, Z. Rudnick and P. Sarnak for useful discussions. The first author acknowledges support from the European Research Council through grant no. 617129.

\section{Review on Cheeger's and Sobolev's constant} \label{s:review}

In this section, we briefly review Cheeger's inequality and the various equivalent definitions of the Cheeger constant. Background on these issues can be found in Chavel's books \cite{chavel,chaveliso} for example and Yau's \cite{yau1975}.

\subsection{Variational characterization of Cheeger's constant}
The Cheeger constant $h(M)$ was defined in $(\ref{cheegdef})$ as an isoperimetric constant. It is possible however to give a variational characterization of $h(M)$ in terms of the \emph{Sobolev constant} $s(M)$ defined by

\begin{equation}\label{eq:sobolev} s(M):=\inf_{f} \frac{||\nabla f||_1}{\inf_\alpha ||f-\alpha||_1}\end{equation}
where the infimums range over Lipschitz functions $f$ on $M$ and real numbers $\alpha$ respectively. It is well-known that
\begin{equation}\label{sobolev}
h(M)=s(M)
\end{equation}

The remainder of this subsection is devoted to a proof of this fact, which we give for the reader's convenience. First observe that the infimum in the right hand side of~\eqref{eq:sobolev} is attained for the number $\alpha_0$ such that $vol (f\leq \alpha_0) = vol (f \geq \alpha_0)$ (for a generic function this number exists and is unique). To see this, it suffices to prove $|| f - \alpha ||_1 \geq || f -\alpha_0 ||_1$ for every $\alpha \geq \alpha_0$ (then the case $\alpha \leq \alpha_0$ is handled by taking $-f$). But decomposing
\[ || f-\alpha ||_1 - ||f-\alpha_0||_1  = \int _{f\geq \alpha} (f-\alpha)  -\int _{f\leq \alpha} (f-\alpha) - \int _{f\geq \alpha_0} (f-\alpha_0) + \int _{f\leq \alpha_0} f-\alpha_0  \]
and using $vol(f\geq \alpha_0) = vol (f \leq \alpha_0)$ gives
\[ || f-\alpha ||_1 - ||f-\alpha_0||_1 = \alpha \big( vol (f\leq \alpha) - vol (f\geq \alpha) \big) - 2 \int _{\alpha_0 \leq f \leq \alpha} f ,  \]
and
\[ || f-\alpha ||_1 - ||f-\alpha_0||_1 \geq \alpha \big( vol (f\leq \alpha) - vol (f\geq \alpha) \big) - 2 \alpha \big( vol (f\leq \alpha ) - vol (f\geq \alpha_0) \big) ,\]
which leads to the desired $||f -\alpha ||_1 - ||f-\alpha_0||_1 \geq 0$ by using $vol (f\geq \alpha_0) = vol (M ) /2$.

We now pass to the proof of $(\ref{sobolev})$. Let $\Omega$ be a domain with smooth boundary of volume bounded by $vol(M)/2$, and $\Omega' = M \setminus \Omega$. For every $\varepsilon>0$ small enough, consider the function $f_\varepsilon : M \rightarrow \mathbb R$ defined by $f_\varepsilon (p) = \pm 1$ if $d(p,\partial \Omega)\geq \varepsilon$, and by $f_\varepsilon (p) = \pm d(p,\partial \Omega) /\varepsilon$ if $d(p,\partial \Omega)\leq \varepsilon$; the sign is positive if $p$ belongs to $\Omega$ and negative otherwise. The function $f_\varepsilon$ is Lipschitz and
\[  || \nabla f_\varepsilon ||_1 = 2 vol (\partial \Omega) + O(\varepsilon).\]
On the other hand, for every $\varepsilon>0$ small enough, we have
\[  \inf_\alpha || f_\varepsilon -\alpha ||_1 \geq 2 vol (\Omega) + O(\varepsilon), \]
so that by taking the limit as $\varepsilon $ tends to $0$,
\[  vol (\partial \Omega ) \geq s(M) \cdot vol (\Omega).\]
This being valid for every $\Omega$, we get $ s(M) \leq h(M)$. To prove the reverse inequality, let $f$ be a generic smooth function, and $\alpha_0$ be the number so that $vol (f\leq \alpha_0) = vol (f\geq \alpha_0)$. The co-area formula gives
\begin{equation}\label{coarea}  \int _{f\geq \alpha_0} || \nabla f || dv  = \int _{\alpha_0} ^\infty  vol (f= t) dt  \end{equation}
and by applying the definition of Cheeger's constant to the domains $\{ f \geq t \}$ for $t\geq \alpha_0$, we get
\[ \int _{f\geq \alpha_0} || \nabla f || dv \geq  h(M) \int _{\alpha_0} ^\infty vol (f\geq t) dt= h(M) \int _{f\geq \alpha_0} |f-\alpha_0| dv .\]
The same argument shows that
\[ \int _{f\leq \alpha_0} || \nabla f || dv \geq  h(M) \int _{f\leq \alpha_0} |f-\alpha_0| dv ,\]
hence
\[  ||\nabla f ||_1 \geq h(M) \inf _\alpha ||f-\alpha||_1.\]
This being valid for every generic smooth function, an approximation argument shows that the same remains true for all Lipschitz functions, and we get $h(M) \leq s(M)$.

\subsection{The Cheeger inequality}
The first non zero eigenvalue of the Laplacian on a compact Riemannian manifold $M$ is defined via the well-known characterization by Rayleigh quotients, namely:

\begin{equation}\label{rayleigh-def}
\lambda_1(M) = \inf_{f} \frac{||\nabla f||^2_2}{ ||f||^2_2} ,
\end{equation}
where the infimum is taken over all Lipschitz functions on $M$ with zero average with respect to the Riemannian volume on $M$.

Cheeger's inequality claims that
$$\lambda_1(M) \geq \frac{1}{4}h(M)^2.$$
Let us recall its proof.  It is again an application of the co-area formula. Let $f: M \rightarrow \mathbb R$ be an eigenvector associated to $\lambda_1(M)$. Taking $-f$ instead if necessary, we may assume that the domain $\Omega = \{ f>0 \}$ has volume $\leq vol (M)/2$. Set $ v := f  ^2_{|\Omega}$. This is a Lipschitz function and $\nabla v = 2 f \nabla f$ on $\Omega$. By Cauchy-Schwarz
\[      ||\nabla v ||_{L^1 (\Omega)}  \leq 2 || f ||_ {L^2 (\Omega)} \cdot || \nabla f ||_{L^2 (\Omega)} \]
and because $f$ vanishes on $\partial \Omega$, Green's formula gives
\begin{equation}\label{upper}  ||\nabla v ||_{L^1(M)} \leq  2 \sqrt{\lambda_1(M)} ||f ||^2_{L^2(\Omega)} .\end{equation}
On the other hand, the co-area formula $(\ref{coarea})$ gives
\[  || \nabla v ||_{L^1(M)} \geq h(M)  || v ||_{L^1(\Omega)} = h(M) ||f ||^2_{L^2(\Omega)} , \]
which, together with~\eqref{upper}, implies Cheeger's inequality.

\section{A strengthened Cheeger inequality for $2$-covers}\label{2cover-sec}
\subsection{Strategy} In this section, we prove Theorem \ref{main}. In a subsequent paper \cite{breuillard-deroin2} we will prove an analogous result for arbitrary finite covers (with the constant $1/4$ replaced by a positive constant $c_d$ depending only on the degree $d$ of the cover). The argument in the general case, where $d\geq 3$ and the cover $M'$ of $M$ may not be Galois, is more involved, in particular if one wishes to find a good constant $c_d$. However the case $d=2$ gives an idea of how to proceed: if $\lambda_1(M')<\lambda_1(M)$ then every real valued eigenfunction $f$ for $\lambda_1(M')$ on $M'$ will be odd and if $\lambda_1(M')$ is very small compared to $\lambda_1(M)$, then $|f|$ will be almost constant on $M'$. So the variation of $f$ will be concentrated around the nodal set and this will then give a upper bound on the Sobolev constant $s(M')$, which yields the desired inequality. We now pass to the details.

\subsection{Proof of Theorem \ref{main}}
We will show that

\begin{equation}\label{2covers}
\lambda_1(M') \geq \min\{\lambda_1(M),\frac{\sqrt{\lambda_1(M)}(1-2\frac{\lambda_1(M')}{\lambda_1(M)})}{2}s(M')\},
\end{equation}
from which Theorem \ref{main} follows easily with the constant $c_2=\frac{1}{4}$, by considering separately the cases $\frac{\lambda_1(M')}{\lambda_1(M)}\leq \frac{1}{4}$ and $\frac{\lambda_1(M')}{\lambda_1(M)}\geq \frac{1}{4}$ and applying $(\ref{sobolev})$ and $(\ref{buser})$.

There is an order $2$ symmetry on $M'$, and $L^2(M')$ splits as an orthogonal direct sum of odd and even functions. The symmetry commutes with the Laplacian, so this decomposition into odd and even functions respects the eigenspace decomposition of the Laplacian. Moreover even functions descend to functions on $L^2(M)$ (they are precisely the lifts to $M'$ of functions on $L^2(M)$). In particular every smooth even function $h$ with zero average on $M'$ satisfies

\begin{equation}\label{even}
\frac{||\nabla h||_2^2}{||h||_2^2} \geq \lambda_1(M) .
\end{equation}
In particular $\lambda_1(M')\leq \lambda_1(M)$ and we may thus assume that $\lambda_1(M')< \lambda_1(M)$. Let $f$ be an eigenfunction of $\Delta$ on $M'$ corresponding to $\lambda_1(M')$. Then clearly $f$ is odd. It is well-known that the nodal set $\{x \in M, f(x)=0\}$ is a closed set of empty interior and zero measure (see e.g. \cite{hardt-simon}).  It follows that $\{x \in M', f(x) \geq 0\}$ is a fundamental domain for the order $2$ symmetry and has measure $\frac{1}{2}vol(M')=vol(M)$, where $vol(M)$ denotes the Riemannian volume of $M$.

The absolute value $|f|$ is an even function on $M'$. It is also Lipschitz and we may write $|f|=\alpha + u$, where $\alpha \in \R$ and $u$ is an even function with zero average on $M'$. We have

\begin{equation}\label{l2dec}
||f||_2^2 = \alpha^2vol(M') + ||u||_2^2 ,
\end{equation}
and

$$||\nabla u ||_2^2 = ||\nabla |f|||_2^2 = ||\nabla f||_2^2 = \lambda_1(M') ||f||_2^2,$$
since $|\nabla |f||= | \nabla f|$ almost everywhere (see e.g. \cite[Lemma VIII.3.2]{chavel}). Applying $(\ref{even})$ we get:
\begin{equation}\label{normbound}
||u||_2^2 \leq \frac{1}{\lambda_1(M)}||\nabla u||_2^2 = \frac{1}{\lambda_1(M)}||\nabla f||_2^2 =   \frac{\lambda_1(M')}{\lambda_1(M)}||f||_2^2 .
\end{equation}
Now set $v$ to be the function on $M'$ defined as $u^2$ on $\{f \geq 0\}$ and $2\alpha^2 - u^2$ on $\{f \leq 0\}$. Note that $v=\alpha^2$ on the nodal set  $\{f=0\}$ and that $v$ is Lipschitz on $M'$ and $\nabla v = \nabla u^2$. Applying  the Cauchy-Schwarz inequality $\nabla u^2 = 2u \nabla u$, the multiplication of the last two inequalities yields:

\begin{equation}
||\nabla v||_1 = ||\nabla u^2||_1 \leq 2||u||_2\cdot ||\nabla u||_2 \leq 2 \frac{\lambda_1(M')}{\sqrt{\lambda_1(M)}} ||f||_2^2 .
\end{equation}

Now consider an arbitrary real number $m \in \R$. We have:

$$||v - m||_1 = \frac{1}{2}(||u^2 - m||_1 + ||2\alpha^2 -u^2 -m||_1) \geq \alpha^2 |M'|- ||u||_2^2$$
by the triangle inequality.

However by $(\ref{l2dec})$ and $(\ref{normbound})$

$$\alpha^2|M'|-||u||_2^2= ||f||_2^2 - 2||u||_2^2 \geq ||f||_2^2(1- 2\frac{\lambda_1(M')}{\lambda_1(M)}). $$

We thus obtain the desired bound:

$$h(M') = s(M') \leq \frac{||\nabla v ||_1}{\inf_m||v -m||_1} \leq 2\frac{\lambda_1(M')}{\sqrt{\lambda_1(M)}(1-2\frac{\lambda_1(M')}{\lambda_1(M)})}.$$

\section{Small Salem numbers and $\lambda_1$}\label{salem-sec}
In this section we recall the proof of Sury's observation $(\ref{kleineq})$ from the introduction, then proceed to prove Theorem \ref{specriterion}.

\subsection{Salem numbers and length of geodesics}

Let us first recall the correspondence between Salem numbers and the length of geodesics as described in \cite{sury,ghate-hironaka,maclachlan-reid}. If $\gamma$ is a hyperbolic element in a Fuschian group $\Gamma$, then the axis of $\gamma$ in $\HH^2$ descends to a closed geodesic on $\HH^2/\Gamma$ of length $2|\log(\lambda_\gamma)|$, where $\lambda_\gamma$ is an eigenvalue of $\gamma$ in $PSL_2(\R)$. If $\Gamma$ is arithmetic, then the field generated by the set of traces $tr(x^2)$, $x \in \Gamma$ is a totally real number field $k_0$ and this set is sent to a bounded subset of $\R$ by any non identity Galois embedding of $k_0$ in $\R$ (see Takeuchi \cite[Theorem 1]{takeuchi}). Considering powers $\gamma^n$, $n$ large, this shows that $\lambda_\gamma^2$ has all its Galois conjugates on the unit circle, except for $\lambda_\gamma^{-2}$, hence is a Salem number.

\begin{rmm}\label{surface2salem} The length of any closed geodesic in an arithmetic $2$-orbifold is equal to $\log \tau$ for some Salem number $\tau>1$.
\end{rmm}

Conversely, suppose $\tau>1$ is a Salem number. Then $k_0:=\Q(\tau+\tau^{-1})$ is a totally real number field, embedded in $\R$. We can build an arithmetic (in fact congruence) Fuchsian group containing a hyperbolic element with eigenvalue $\tau$ as follows. Quaternion algebras over $k_0$ are determined up to isomorphism by the set $Ram(A)$ of places of $k_0$ where they ramify. This set is always of even cardinality and to any such set corresponds a unique quaternion algebra. If $[k_0:\Q]$ is odd, we take $Ram(A)$ to be the set of all archimedean places of $k_0$ except its defining real place. If $[k_0:\Q]$ is even, we will add just one finite place, a large prime $\P$ of $k_0$ to be specified below. Let $A$ be the corresponding quaternion algebra of $k_0$, and $A^1$ be the subgroup of elements of reduced norm $1$.

We want to find $\P$ in such a way that $\Q(\tau)$ embeds in $A$ as a $k_0$-algebra. It is a well-known property of quaternion algebras (see e.g. \cite[Theorem 7.3.3]{maclachlan-reid}) that a quadratic extension $k$ of $k_0$ embeds in $A$ iff each place of $k_0$ in $Ram(A)$ is inert, i.e. there is only one place of $k$ above each place of $k_0$ in $Ram(A)$. Since $\tau$ is Salem, $k_0(\tau)|k_0$ has just one (complex) place above each real place other than the original one. By Tchebotarev's density theorem, there are infinitely many inert prime ideals $\P$ in $k_0$. Choose one to define $Ram(A)$ when $[k_0:\Q]$ is even. We have thus constructed $A$ in such a way that $k:=\Q(\tau)$ embeds in $A$. If $(a,b)$ is a $k$-basis of $A$, then $I=\mathcal{O}_ka+\mathcal{O}_kb$ is an integral ideal of $A$ and the set $\{\alpha \in A; \alpha I \subset I\}$ is an order containing $\tau$. Let $\mathcal{O}$ be a maximal order containing it. Then $\mathcal{O} \cap A^1$ is a discrete subgroup of finite covolume in $A^1(\R) \simeq \SL_2(\R)$. It contains $\tau$, because $N_{\Q(\tau)|k_0}(\tau)=1$ and as $\Q(\tau)$ embeds in $A$ the norm $N_{\Q(\tau)|k_0}$ coincides with the reduced norm of $A$. Denote by $\Gamma_\mathcal{O}$ its image in $PSL_2(\R)$. Since $\tau>1$, the associated element in $\Gamma_\mathcal{O}$ is hyperbolic. We conclude:

\begin{rmm}\label{salem2surface} For every Salem number $\tau>1$, there is a cocompact arithmetic Fuchsian group $\Gamma_{\mathcal{O}}$ associated to a maximal order $\mathcal{O}$ in a quaternion algebra defined over $\Q(\tau+\tau^{-1})$, such that $2\log \tau$ is the length of a closed geodesic in the $2$-orbifold $\HH^2/\Gamma$.
\end{rmm}

In general $\mathcal{O}$ may have torsion, but, as shown to us by Alan Reid,  it is possible to choose the quaternion algebra over $k_0$  in such a way that $A^1$ contains no non trivial torsion element. Indeed if $A^1$ contains an element of order $n$, then the number field $k_0(e^{\frac{2i\pi}{n}})$ is a quadratic extension of $k_0$ and it embeds as a $k_0$ subalgebra of $A$. However, as before, this is the case if and only if every place in $Ram(A)$ is inert with respect to this quadratic extension (see \cite[Theorem 12.5.4]{maclachlan-reid}). Again by Tchebotarev's density theorem, there are infinitely many split primes $P$ in $k_0$. So we can find such a prime ideal $\P_n$. On the other hand, since $\cos(2\pi/n) \in k_0$, the degree of $e^{\frac{2i\pi}{n}}$ over $\Q$ is bounded, and hence there is an upper bound on the order $n$ of torsion elements. We may thus modify $Ram(A)$ by adding finitely many prime ideals of the form $\P_n$. This guarantees that $A^1(k_0)$ is torsion free. We conclude:

\begin{rmm}\label{salem2smoothsurface} For every Salem number $\tau>1$, there is a torsion-free cocompact arithmetic congruence Fuchsian group $\Gamma$, such that $2\log \tau$ is the length of a closed geodesic in the hyperbolic surface $\HH^2/\Gamma$.
 \end{rmm}
In particular:
\begin{rmm}
Salem's conjecture is equivalent to a uniform lower bound on the length of closed geodesics in congruence arithmetic hyperbolic surfaces (without conical points).
\end{rmm}

\subsection{Short geodesics and $\lambda_1$}
Here we prove Theorem $(\ref{specriterion})$ from the introduction. Assume that Salem's conjecture holds. By $(\ref{surface2salem})$ there is a uniform lower bound, say $c_0>0$ on the length of closed geodesics in arithmetic $2$-orbifolds. Now suppose $\Sigma=\HH/\Gamma$ is a congruence arithmetic $2$-orbifold. Its $\lambda_1$ (for orbifolds it is defined as the infimum of the Rayleigh quotients $(\ref{rayleigh-def})$ over all $\Gamma$-invariant Lipchitz functions on $\HH^2$) is bounded away from zero uniformly, thanks to the known approximations towards the Ramanujan conjecture for $SL_2$ over an arbitrary number field, which follow from the works of Gelbart-Jacquet \cite{gelbart-jacquet} by applying the Jacquet-Langlands correspondence, and extend to compact arithmetic surfaces the work of Selberg \cite{selberg} for non compact ones. See Vigneras \cite{vigneras}, where the following is explained:

\begin{rmm}\label{selberg} Let $\Sigma=\HH^2/\Gamma$ be a congruence arithmetic $2$-orbifold. Then
$$\lambda_1(\Sigma) > \frac{3}{16}.$$
\end{rmm}

Now pick a $2$-cover $\widetilde{\Sigma}$ of $\Sigma$. We wish to apply Theorem \ref{main} to $\widetilde{\Sigma}$. In order to do so, note that the proof given in Section \ref{2cover-sec} extends verbatim to the case of $2$-orbifolds. We conclude that

\begin{equation}\label{ramabound}\lambda_1(\widetilde{\Sigma}) \geq \frac{1}{4} \sqrt{\lambda_1(\Sigma)} \cdot h(\widetilde{\Sigma}) \geq \frac{\sqrt{3}}{16} h(\widetilde{\Sigma}),
\end{equation}
where $h(\widetilde{\Sigma})$ is the Cheeger constant:
$$h(\widetilde{\Sigma})= \inf \frac{\ell(\mathcal{C})}{\min \{area(B), area(\widetilde{\Sigma} \setminus B)\}},$$
where $\mathcal{C}$ ranges over closed smooth curves separating $\widetilde{\Sigma}$ into two disjoint connected open sets $B$ and $\widetilde{\Sigma} \setminus B$ (see \cite[\S 8.3]{buser}). We claim that there is a constant $c_1>0$ independent of $\widetilde{\Sigma}$ such that
\begin{equation}\label{hbound}
h(\widetilde{\Sigma}) \geq \frac{c_1}{area(\widetilde{\Sigma})}.
\end{equation}
Combined with $(\ref{ramabound})$ this yields the lower bound in Theorem $(\ref{specriterion})$ as desired.
In order to prove the claim, let $\mathcal{C}$ be a smooth closed curve almost realizing the infimum in the definition of $h(\widetilde{\Sigma})$. If $\mathcal{C}$ is not null homotopic, then its length is always at least equal to the length of a closed geodesic: indeed $\mathcal{C}$ lifts to $\HH^2$ to a curve between two distinct points say $x_0$ and $\gamma x_0$ for some \emph{hyperbolic} $\gamma \in \Gamma$, hence its length $\ell(\mathcal{C})$ must be at least the minimal displacement of $\gamma$, i.e. $\inf_{x \in \HH^2} d(x,\gamma x)$,  hence at least $c_0$ by $(\ref{surface2salem})$ and our assumption that Salem's conjecture holds. Hence $(\ref{hbound})$ holds in this case. On the other hand, if $\mathcal{C}$ is null homotopic, then it lifts to a curve in $\HH^2$ between two points say $x_0$ and $\gamma x_0$, with $\gamma$ an elliptic (possibly trivial) isometry of $\HH^2$. Let $T$ be the triangle with base point the elliptic fixed point $p$,  two geodesic sides connecting $p$ to $x_0$ and $\gamma x_0$ and a third side equal to $\mathcal{C}$ (if $x_0=\gamma x_0$ the lift of $\mathcal{C}$ bounds a disc: then pick $p$ inside the disc and $T$ becomes degenerate). The computation in \cite[8.4.1]{buser} (see also \cite[Prop.3]{yau}) shows that $\ell(\mathcal{C}) \geq area(T)$. This yields $h(\widetilde{\Sigma}) \geq 1$, which in particular gives $(\ref{hbound})$.

\bigskip

We now pass to the converse in the proof of Theorem $(\ref{specriterion})$, namely we assume the existence of the lower bound $\lambda_1(\widetilde{\Sigma})$ for all $2$-covers $\widetilde{\Sigma}$ and deduce Salem's conjecture. By $(\ref{salem2surface})$ Salem's conjecture follows if we can prove a uniform lower bound on the length of closed geodesics of all congruence arithmetic $2$-orbifolds. In fact  it is enough to consider torsion free groups of the form $\Gamma_{\mathcal{O}}$ defined in $(\ref{salem2surface})$. 

Let $\Sigma:=\HH^2/\Gamma_{\mathcal{O}}$, where $\Gamma_{\mathcal{O}}$ was defined above. Let $\mathcal{C}$ be a closed geodesic in $\Sigma$ corresponding to a hyperbolic element $\gamma \in \Gamma_{\mathcal{O}}$ with axis $a_\gamma$. We say that $\mathcal{C}$ is \emph{simple} if there is no $\delta \in \Gamma$ such that $\delta(a_\gamma)$ self-intersects $a_\gamma$ non trivially (i.e. we do not allow $\delta(a_\gamma) = a_\gamma$, which may happen if $\delta$ is an elliptic element of order $2$, or if some power of $\delta$ coincides with a power of $\gamma$). It is shown in Beardon's book \cite[Theorem 11.6.8(1)]{beardon} that this situation, the length of $\mathcal{C}$ is at least $0.2$. Hence we may assume that $\mathcal{C}$ is a simple closed geodesic.

We cut along this geodesic. If this operation separates $\Sigma$ into two disjoint connected components, then $h(\Sigma) \leq \frac{\ell(\mathcal{C})}{area(B)}$, where $B$ is the component of smaller area. Note that $area(B)$ is bounded away from zero: indeed pasting together two copies of $B$ along $\mathcal{C}$ gives rise to a compact $2$-orbifold, but those have a minimal co-covolume by the volume formula for fuchsian groups (the $(2,3,7)$ triangle group realizes the minimum).  Combined with Buser's inequality (\ref{buser}), and the spectral bound $(\ref{ramabound})$, this gives an explicit uniform lower bound on $\ell(\mathcal{C})$.

Now if cutting along $\mathcal{C}$ does not separate $\Sigma$, we may paste together two copies of $\Sigma$ by gluing them along the cut and obtain in this way another compact connected $2$-orbifold $\widetilde{\Sigma}$ (see \cite[Chapter 11]{ghysharpe} for the uniformization of orbifolds), which is a $2$-cover of $\Sigma$. Its Cheeger constant satisfies
$$h(\widetilde{\Sigma}) \leq \frac{\ell(\mathcal{C})}{area(\Sigma)}.$$
Hence by Buser's inequality $(\ref{buser})$
$$\lambda_1(\widetilde{\Sigma}) \leq \frac{\ell(\mathcal{C})}{area(\Sigma)}.$$
This combined with our assumption on the lower bound for $\lambda_1(\widetilde{\Sigma})$ gives the desired uniform lower bound on $\mathcal{C}$ and this ends the proof of Theorem $(\ref{specriterion})$.

\subsection{Hyperbolic $n$-manifolds, $n \geq 3$}

In this paragraph, we deal with $(\ref{kleineq})$, which claims that in dimension $3$ and higher the uniform lower bound in $1/volume$ for the $\lambda_1$ of $2$-covers of compact congruence arithmetic $n$-manifolds holds true without assumptions. As we now explain, this follows from the combination of the Burger-Sarnak uniform lower bound on $\lambda_1$ for congruence arithmetic hyperbolic $n$-manifolds, our Theorem \ref{main} for $2$-covers and a result of Schoen \cite{schoen} pertaining to lower bounds for the Cheeger constant of negatively curved manifolds.

Schoen \cite{schoen} proved that if $M$ is a compact hyperbolic $n$-manifold with $n \geq 3$, then its Cheeger constant $h(M)$ is bounded below by $c_n/vol(M)$, where $c_n$ is a positive constant depending only on the dimension $n$. In fact Schoen stated his result in terms of a lower bound for $\lambda_1(M)$ of the form $c_n/vol(M)^2$, but his proof deals with the Cheeger constant and only at the end Cheeger's inequality is applied to give the lower bound on $\lambda_1$. Combining this with our Theorem \ref{main}, we see that if $M'$ is any $2$-cover of $M$, then
\begin{equation}\label{sch}\lambda_1(M') \geq (c_n/4)\frac{\sqrt{\lambda_1(M)}}{vol(M')}.\end{equation}

Burger-Sarnak \cite{burger-sarnak} proved a restriction principle for automorphic representations of algebraic groups, which allows to deduce a spectral gap for the automorphic representation of the ambient group from a similar spectral gap for automorphic representations of an algebraic subgroup. Using the Gelbart-Jacquet \cite{gelbart-jacquet} spectral gap for automorphic representations of $SL_2$ over an arbitrary number field already mentioned above, they could apply their restriction principle to deduce \cite[Corollary 1.3]{burger-sarnak} the following spectral gap for the $\lambda_1$ of compact hyperbolic congruence $n$-manifolds:
\begin{equation}\label{bs}\lambda_1(M) \geq \frac{2n-3}{4}.
\end{equation}
Actually they only stated this bound for hyperbolic manifolds derived from a quadratic form. However their method works in the general case as well (all associated $K$-forms of $SO(n,1)$ admit an embedded isogenous copy of $SL_2$). We thank Nicolas Bergeron for this observation. More recently Bergeron and Clozel, making use of J. Arthur's work, improved considerably the above spectral gap \cite{bergeron-clozel}. For us the Burger-Sarnak bound is enough.
Now $(\ref{bs})$ together with $(\ref{sch})$ yields the desired lower bound in $(\ref{kleineq})$, namely:
$$\lambda_1(M') \geq c'_n/vol(M'),$$ for some $c'_n>0$ depending on $n$ only.
As already mentioned in the introduction, one can build for each $n \geq 2$ infinite families of $2$-covers of congruence hyperbolic $n$-manifolds such that $\lambda_1(M') \leq C_n / vol(M')$, showing that lower bound in $(\ref{kleineq})$ is sharp up to multiplicative constants.

\bibliographystyle{abbrv}

\end{document}